\documentclass[10pt]{amsart}       %{proc-l}
\usepackage{amssymb}
\newtheorem{theorem}{Theorem}[section]

\newtheorem{proposition}{Proposition}[section]   % not in PAMS

\theoremstyle{definition}

\theoremstyle{remark}

\numberwithin{equation}{section}

\newcommand{\on}{\operatorname}
\begin{document}
%\magstep2
%Topmatter 

%One author
\title[On a Combinatorial Identity]
    {On a Combinatorial Identity}% \\%<title line 2>

\author{Jintai Ding}
\address{Ding: Department of Mathematical Sciences,
  University of Cincinnati, 
   Cincinnati, OH 45221}

\email{ding@math.uc.edu}
\author{Naihuan Jing}
\address{Jing: Mathematical
Sciences Research Institute,
1000 Centennial Drive, Berkeley, CA
94720-5070}
\address{Department of Mathematics,
North Carolina State Univer\-sity, 
   Ra\-leigh, NC 27695-8205}

\email{jing@math.ncsu.edu}
\thanks{Jing gratefully acknowledges support from
NSA grant MDA904-97-1-0062 and Mathematical Sciences Research Insititute
through NSF grant DMS-9701755}
%\date{}

%End topmatter

\begin{abstract}
Recently the second named
author discovered a combinatorial
identity in the context of vertex representations of quantum
Kac-Moody algebras. We give a direct and
elementary proof of this identity.
Our method is to show a related
identity of distributions.
\end{abstract}

\maketitle

%\centerline{To appear}

\section{Introduction}
\label{S:intro}

In \cite{J2} one of us  proved the following combinatorial
identity
in the context of quantum Kac-Moody algebras.
For any integer
$m\geq 0$ we have

\begin{multline} \label{ident1}
\sum_{\sigma\in
S_{m+1}}\sum_{r=0}^{m+1}\sigma.
\begin{bmatrix}
m+1\\r\end{bmatrix}(w-q^{m}z_1)\cdots (w-q^{m}z_r)\\
\cdot
(z_{r+1}-q^{m}w)\cdots
(z_m-q^{m}w)\prod_{i<j}
\frac{z_i-q^{2}z_j}{z_i-z_j}=0,
\end{multline}
where
 the symmetric group $S_{m+1}$ acts on the variables $z_i$ by permuting
the
indices while fixing $q$ and $w$. Here the $q$-Gaussian number
$\begin{bmatrix} m+1\\r\end{bmatrix}$ is defined as
follows.
\begin{equation*}
\begin{bmatrix}
m+1\\r\end{bmatrix}=\frac{[m+1]!}{[r]![m+1-r]!},
\qquad [n]!=\prod_{1\leq
i\leq n}[i], \qquad [i]=\frac{q^i-q^{-i}}{q-q^{-1}},
\end{equation*}
and
$[0]=1$.
Observe that this identity is valid in the ring of polynomials
in
$z_1, \cdots, z_{m+1}$ over $\mathbb Z[q, q^{-1}]$.

From combinatorial
viewpoint the identity (\ref{ident1}) 
is equivalent to $m+2$ identities of linear relations among
Hall-Littlewood polynomials \cite{M}
associated to certain tuples (not necessary partitions).
The first and the last of these linear relations are actually equivalent 
to the well-known $q$-binomial identity:
\begin{equation*}
\sum_{r=0}^{m+1}(-1)^r
\begin{bmatrix} m+1\\r\end{bmatrix}=0.
\end{equation*}

This identity was proved in \cite{J2} by interpreting it
as the Serre relation with the
help of vertex 
representation of the quantum Kac-Moody algebra. 
The special cases of $m=1, 2, 3$ were known earlier and 
used in the level one
vertex representations of quantum affine algebras \cite{FJ, J1, J2}.

Later Tarasov
\cite{T} generalized the combinatorial 
identity in the
context of the elliptic quantum algebra
and he proved a more general
identity using a different method. 
In \cite{E} Enriquez proposed a
distribution ``generalization''
of the identity using the quantum shuffle
algebra and he verified
the distribution identity for $m=1, 2$. These cover
all 
the quantum shuttle algebra associated the quantum affine
algebras
except $U_q(G_2^{(2)})$, which corresponds to the distribution
identity
of $m=3$.

The purpose of this paper is to give an elementary
(combinatorial)
proof of the
identity (\ref{ident1}) for arbitrary $m$. The key of  
our
proof is  that the distribution identity is equivalent to
the
combinatorial identity (\ref{ident1}).

To state our result we let the
formal delta-function $\delta(z, w)$ be the following infinite
series:
\begin{equation*}
\delta(z, w)=\sum_{n\in \mathbb
Z}z^{i-1}w^i
\end{equation*}
We also make the convention that the rational
function
$\frac 1{z-w}$ represents the formal series in the direction
$|z|>|w|$: $\sum_{n=0}^{\infty} z^{-1}(\frac wz)^n$. Then we
have
\begin{equation}
\delta(z, w)=\frac 1{z-w}+\frac 1{w-z}
\end{equation}

Using (1.1), we will prove the following theorem.
\begin{theorem}
\label{thm}
For any non-negative integer number $m$ we have:
\begin{align}
\label{ident2} \nonumber
  & {\sum}_{\sigma\in S_{m+1}}
\sigma.\sum_{k=0}^{m+1} \bmatrix m+1 \\ k
  \endbmatrix \frac 1{(q^{-m}z_1
- w)\cdots (q^{-m}z_k -
      w)}\\
&\qquad\qquad
\cdot\frac1{(q^{-m}w-z_{k+1})\cdots (q^{-m}w-z_{m+1})}
{\prod_{i<j}\frac{z_i
      - z_j}{q^2 z_i - z_j}}    
\\ & = q^{m-1}
{\sum}_{\sigma\in S_{m+1}} \sigma.\delta(w,q^{-m}z_1)\delta(z_1,q^2 z_2)
\cdots \delta(z_{m},q^2 z_{m+1}) \nonumber
\end{align}
where the symmetric
group $S_{m+1}$ acts on the indices of the $z_i$.
\end{theorem}
The above
identity is of great importance in representation theory, 
because it
provides an analytical method to describe the Serre relations
for quantum
affine algebras \cite{E}. As an example
the distribution identity for $m=3$
completes 
Eneriquez's argument for the exceptional case of
$U_q(G_2^{(1)})$.

We will see that each step in our proof of theorem
\ref{thm}
is reversible, which means that
we also obtain a new and elementary proof of
(\ref{ident1}).

\section{Delta-function identity} \label{S:delta}
The
delta function plays an important role in the theory of vertex 
operator
algebras and is characterized by the following important 
property. For any
formal distribution
$f(z)\in \mathbb C[[z, z^{-1}]]$ one has
that
\begin{equation*}
f(z)\delta(z, w)= f(w)\delta(z, w)
\end{equation*}

To simplify the presentation we will use 
${\on{Sym}_{{S_{m+1}}}}$ to the
denote the symmetrization operator $\sum_{\sigma\in S_{m+1}}\sigma$ on the
ring of formal series in $z_1, z_2, \cdots,
z_{m+1}$. Let $\omega_i$
($i\geq 2$) 
be the transposition $(1i)$, 
and $S_{m}$ be the symmetric
group on letters $2, 3, \cdots, m+1$. Then we have the coset decomposition
of $S_{m+1}$:
$S_{m+1}=\cup_{i=1}^{m+1} \omega_iS_{m}$, where
$\omega_1=1$.

On the other hand, we need to do some preparation to
interpret the formal power series in a slightly different way. 
Recall
that the function $\frac1{z-w}$ represents a formal power
series in the
region $|z|>|w|$, which we say for simplicity that
$\frac1{z-w}$ lives in
the region $|z|>|w|$.
We now proceed the argument in several
steps.

%\medskip
{\bf Let us now fix the following condition on $q$}:
\begin{equation} 
{  |q|>1 } 
\end{equation} 

Let $0<i<j$. We know that
the function   
$\frac 1{q^{2}z_i-z_j}$ lives in the region
$|q^{2}z_i|>|z_j|$, which is equivalent to the
region
$|z_i|>|z_j||q^{-2}|$. This region contains the subset of
$|z_i|>|z_j|$. 
Let us denote the region $|z_i|>|z_j||q^{-2}|$  by
$R_{ij}$. 

Note that 
the function   
$\frac 1{q^{2}z_j-z_i}$ lives in the
region 
$|q^{2}z_j|>|z_i|$, which is equivalent to the region
$|z_j|>|z_i||q^{-2}|$.
We denote this region by $R_{ji}$. Clearly the subset  
$|z_j|>|z_i|$ is
contained in $R_{ji}$.

On the other hand $R_{ij} \bigcap R_{ji}$ is not
empty: 
\begin{equation}
R_{ij} \bigcap R_{ji} =\{ z_i, z_j,
|q^2||z_j|>|z_i|>|z_j||q^{-2}| \}.  
\end{equation}

This says that there is a common region $R_{ij} \bigcap R_{ji}$, in which
both 
$\frac 1{q^{2}z_j-z_i}$ and $\frac 1{q^{2}z_i-z_j}$ live.

\medskip
%\vskip .5in 

On the other hand, the function   
$\frac
1{q^{-m}z_i-w}$ lives in the region 
$|q^{-m}z_i|>|w|$, which is equivalent
to the region $|z_i|>|w||q^m|$. 
We denote this region by $R_{iw}.$
The
function   
$\frac 1{q^{-m}w-z_i}$ lives in the region 
$|q^{-m}w|>|z_i|$,
which is equivalent to the region $|w|>|z_i||q^m|$.
We denote this region
by  $R_{wi}.$

Clearly we have a situation different from  what we had
above: 
\begin{equation} 
R_{wi}\bigcap R_{iw}= \emptyset.
\end{equation}

{\bf Step 1}

Let's first change the expansion direction for
the series
$\frac 1{(q^{-m}z_1 - w)\cdots (q^{-m}z_k -
      w)}$ by
using
\begin{equation}
\frac 1{q^{-m}z_i-w}=-\frac 1{w-q^{-m}z_i}+\delta(w,
q^{-m}z_i).
\end{equation}

Then we have 
\begin{align*}
\mbox{LHS}
&=\mbox{LHS of Eqn (1.3)= \label{ident1}}\\
& {\sum}_{\sigma\in S_{m+1}}
\sigma.\sum_{k=0}^{m+1} \bmatrix m+1 \\ k
  \endbmatrix (-1)^k\frac
1{(w-q^{-m}z_1)\cdots (w-q^{-m}z_k 
      )}\\
&\qquad\qquad
\cdot\frac1{(q^{-m}w-z_{k+1})\cdots (q^{-m}w-z_{m+1})}
{\prod_{i<j}\frac{z_i
      - z_j}{q^2 z_i - z_j}} +  \\
&\underset{{S_{{m+1}}}}
{\mbox{\large Sym}}
\sum_{k=0}^{m+1} \bmatrix m+1
\\ k
  \endbmatrix  \sum_{l=1}^k
\frac {\delta(w, q^{-m}z_l)}{(q^{-m}z_1 -
w)\cdots
(w-q^{-m}z_{l-1})}\cdot\\
&\cdot\frac{(-1)^{k-1}}{(w-q^{-m}z_{l+1})
(q^{-m}z_k -
      w)(q^{-m}w-z_{k+1})\cdots
(q^{-m}w-z_{m+1})}
{\prod_{i<j}\frac{z_i
      - z_j}{q^2 z_i -
z_j}}\\
&=X_1+\sum_{n=1}^{m+1}(1,n)\underset{{S_{{m}}}}
{\mbox{\large
Sym}}
\sum_{k=0}^{m+1} \bmatrix m+1 \\ k
  \endbmatrix\sum_{l=1}^k
\frac
{\delta(w, q^{-m}z_1)}{(q^{-m}z_1 - w)\cdots (w-q^{-m}z_{l-1})}\cdots\\
=&
{X_1} +\sum_{n=1}^{m+1}\sum_{k=0}^{m+1}
\bmatrix m+1 \\ k
\endbmatrix
\delta(q^{-m}z_n,
w)\sum_{l=1}^k(1,n)\underset{{S_{{m}}}}
{\mbox{\large
Sym}}(-1)^{k-1}(1,l)\cdot\\
&\hskip 1in \frac 1{(q^{-m}z_1 - w)\cdots
(w-q^{-m}z_{l-1})}\cdots\\
\end{align*}
where we have decomposed the action
of $S_{m+1}$ into the summation
of the cosets of $S_m$. Here $S_m$ is the
symmetric group
of the letters $2, \dots, m$; and by 
${X_1}$ %${\bf X}$
we
mean the  first summation where 
there is no delta functions. 

From the
preparation of expansion region, we see that  
the term $X_1$  lives in the
region:
\begin{equation}
 Y_1= R_{w1}\bigcap R_{w2} \bigcap ...
R_{w(m+1)}
\bigcap \{z_1,...z_{m+1}|\bigcap_{0<i<j<m+1}( R_{i,j}\bigcap R_{ji})\}.
\end{equation}

It follows from (\ref{ident1}) that on $Y_1$ we
have
\begin{equation}
X_1=0.
\end{equation}

Then we have
\begin{align*}
\mbox{LHS} &=
\\
&\sum_{n=1}^{m+1}\sum_{k=0}^{m+1}\bmatrix m+1 \\ k
\endbmatrix
\delta(q^{-m}z_n,
w)\sum_{l=1}^k(1,n)\underset{{S_{{m}}}}
{\mbox{\large
Sym}}(-1)^{k-1}(1,l)\cdot\\
&\hskip 1in \frac 1{(q^{-m}z_1 - w)\cdots
(w-q^{-m}z_{l-1})}\cdots\\
\end{align*}

Now let us take a look at the
first term of the summation. We have

\begin{align*}
&\sum_{k=0}^{m+1}
\bmatrix m+1 \\ k
\endbmatrix
\delta(q^{-m}z_1,
w)\sum_{l=1}^k\underset{{S_{{m}}}}
{\mbox{\large
Sym}}\frac{(-1)^{k-1}}
{(w-q^{-m}z_2)\cdots (w-q^{-m}z_k)}\cdot\\
&\hskip 1in \cdot\frac 1{(q^{-m}w-z_{k+1})\cdots
(q^{-m}w-z_{m+1})}
(1,l){\prod_{i<j}\frac{z_i
      - z_j}{q^2 z_i -
z_j}}\\
&=\delta(q^{-m}z_1, w)\sum_{k=0}^{m+1}
\bmatrix m+1 \\ k
\endbmatrix
\sum_{l=1}^k\underset{{S_{{m}}}}
{\mbox{\large
Sym}}\frac{(-1)^{k-1}q^{m(k-1)}}
{(z_1-z_2)\cdots (z_1-z_k)}\cdot\\
&\hskip 1 in \cdot\frac 1{(q^{-2m}z_1-z_{k+1})\cdots
(q^{-2m}z_1-z_{m+1})}
(1,l){\prod_{i<j}\frac{z_i
      - z_j}{q^2 z_i -
z_j}}\\
&=\delta(q^{-m}z_1, w)\sum_{k=0}^{m+1}
\bmatrix m+1 \\ k
\endbmatrix
\sum_{l=1}^k\underset{{S_{{m}}}}
{\mbox{\large
Sym}}\frac{(-1)^{k-1}q^{m(k-1)}}
{(z_1-z_2)\cdots (z_1-z_k)}\cdot\\
&\hskip
1 in \cdot\frac 1{(q^{-2m}z_1-z_{k+1})\cdots (q^{-2m}z_1-z_{m+1})}
(1,2, \dots, l){\prod_{i<j}\frac{z_i
      - z_j}{q^2 z_i -
z_j}}
\end{align*}
since $(1, 2, \dots, l)(1,2)$ is an element in
$S_{m}=S_m(z_2,\dots, z_{m+1})$.
Then the above is simplified to the
following expression:
\begin{align*}
&\delta(q^{-m}z_1,
w)\sum_{k=0}^{m+1}
\bmatrix m+1 \\ k
\endbmatrix
\sum_{l=1}^k\underset{{S_{{m}}}}
{\mbox{\large
Sym}}\frac{(-1)^{k-1}q^{m(k-1)}(z_1-z_{k+1})\cdots
(z_1-z_{m+1})}
{(q^{-2m}z_1-z_{k+1})\cdots
(q^{-2m}z_1-z_{m+1})}\cdot\\
&\quad\cdot\frac 1{(q^{2}z_2-z_1)\cdots
(q^{2}z_l-z_1)
(q^2z_1-z_{l+1})\cdots (q^2z_1-z_{m+1})}{\prod_{2\leq
i<j\leq m+1}\frac{z_i
      - z_j}{q^2 z_i - z_j}}, 
\end{align*}
which we
will denote by
$\delta(q^{-m}z_1, w) X_{1,w}$. 

It is clear as an analytic
function  
\begin{equation} 
X_{1,w}= \frac 1{2\pi \sqrt {-1}}
\int_{\mbox{w around} q^{-m}z_1} 
-\mbox{LHS } 
d_w.
\end{equation}

Therefore based upon (\ref{ident1}) we have obtained the
following result. 
\begin{proposition}
As an analytic function,

\begin{equation}
X_{1,w}=0
\end{equation}
\end{proposition}

{\bf Step
2}

Now let's change the expansion direction of $z_1$ in
$\frac1{q^2z_i-z_1}$
for $i=2, \dots, l$. Following the same argument as in
{\bf Step 1},
and we see that the 
above expression
becomes:
\begin{align*}
&\delta(q^{-m}z_1, w)\sum_{k=0}^{m+1}
\bmatrix m+1
\\ k \endbmatrix
\sum_{l=1}^k\underset{{S_{{m}}}}
{\mbox{\large
Sym}}\frac{(-1)^{k-1}q^{m(k-1)}(z_1-z_{k+1})\cdots
(z_1-z_{m+1})}
{(q^{-2m}z_1-z_{k+1})\cdots
(q^{-2m}z_1-z_{m+1})}\cdot\\
&\quad\cdot\sum_{l_1=1}^l
\delta(z_1,
q^2z_{l_1})\frac {(-1)^{l-1}}
{(z_1-q^{2}z_2)\cdots
(z_1-q^{2}z_l)
(q^2z_1-z_{l+1})\cdots (q^2z_1-z_{m+1})}\\
&\hskip 2in
{\prod_{2\leq i<j\leq m+1}\frac{z_i
      - z_j}{q^2 z_i -
z_j}}\\
&=\delta(q^{-m}z_1, w)\sum_{k=0}^{m+1}
\bmatrix m+1 \\ k
\endbmatrix
\sum_{l=1}^k\underset{{S_{{m}}}}
{\mbox{\large
Sym}}\frac{(-1)^{k-1}q^{m(k-1)}(z_1-z_{k+1})\cdots
(z_1-z_{m+1})}
{(q^{-2m}z_1-z_{k+1})\cdots
(q^{-2m}z_1-z_{m+1})}\cdot\\
&\quad\cdot\sum_{l_1=1}^l
\delta(z_1,
q^2z_{l_2})\frac {(-1)^{l-1}}
{(z_1-q^{2}z_2)\cdots
(z_1-q^{2}z_l)
(q^2z_1-z_{l+1})\cdots (q^2z_1-z_{m+1})}(2,l_1)\\
&\hskip
2in {\prod_{2\leq i<j\leq m+1}\frac{z_i
      - z_j}{q^2 z_i -
z_j}}\\
\end{align*}

{\bf Step 3} 
After this, we can use the same
argument to change the expansion 
direction of $z_2$, and then
$z_3$ and so
on.  The final term will be  $q^{m-1}\delta(z_m,q^2z_{m+1})$, 
which gives
us a proof of (1.3). 

This gives us the proof under the condition that
$q>1$. 
Using the argument of analytic continuation, 
it is clear that
(1.3) also holds for the case that $q<1$. 
This gives us the proof for our
Theorem.

\medskip
%\vskip .2in 
If we observe carefully the proof
above,  we can see that each step step 
above is reversible. This means
that we 
can start with the right hand side 
of the delta identity and
follow  the proof steps backward to the
beginning.  
After finishing Step 1
in the end,  replace the
delta function by $0$, which actually gives us a
new proof of identity (1.1). 
This is a direct and elementary proof of
identity (1.1).    
\vskip .2in

 \noindent {\it Remark.}  
By replacing
the delta function by $0$, 
we also  obtain an elementary proof of the
combinatorial identity (6.1) of \cite{J2}.

\end{document}